\documentclass[a4paper]{article}
\usepackage{charter}
\usepackage{mathptmx}
\usepackage{amssymb}
\newcommand{\ds}{\displaystyle}
\newcommand{\coco}{\vspace{-.1mm}}

\def\buildrel#1_#2^#3{\mathrel{\mathop{\kern 0pt#1}\limits_{#2}^{#3}}}

\def\cprime{$'$}

\newtheorem{dfn}{Definition}%[section] 
\newtheorem{rmk}{Remark}%[section]
\newtheorem{thm}{Theorem}%[section] 
\newtheorem{cor}{Corollary}%[section]
\newtheorem{prop}{Proposition}%[section] 
\newtheorem{lem}{Lemma}%[section]
\newcommand{\Pf}{{\em Proof}. }
\newcommand{\SPf}{{\em Sketch of proof}. }
\newcommand{\EPf}
{%
\mbox{}%
\nolinebreak%
\hfill%
\rule{2mm}{2mm}%
\medbreak%
\par%
}

\newcommand{\End}{\mbox{$\mathsf{End}$}}
\newcommand{\Hom}{\mbox{$\mathsf{Hom}$}}

\newcommand{\Id}{\mbox{$\mathsf{Id}$}}
\newcommand{\Pol}{\mbox{$\mathsf{Pol}$}}

\begin{document}

\title{Topological Hopf algebras, quantum groups \\
and deformation quantization}

\author{Philippe Bonneau and Daniel Sternheimer\\
\textsl{\small Laboratoire Gevrey de Math\'ematique 
physique, Universit\'e de Bourgogne,}\\ 
\textsl{\small BP 47870, F-21078 Dijon Cedex, France.}\\
\texttt{\small Philippe.Bonneau@u-bourgogne.fr, 
Daniel.Sternheimer@u-bourgogne.fr}}

\date{\small January 5, 2003; revised May 7, 2003}

\maketitle

\begin{abstract}
After a presentation of the context and a brief reminder of deformation 
quantization, we indicate how the introduction of natural topological 
vector space topologies on Hopf algebras associated with Poisson Lie groups, 
Lie bialgebras and their doubles explains their dualities and provides a 
comprehensive framework. Relations with deformation quantization and 
applications to the deformation quantization of symmetric spaces are described.

\end{abstract}

\noindent \textbf{Mathematics Subject Classifications (2000)}: 
Primary 54C40, 14E20; Secondary 46E25, 20C20\\
\textbf{Keywords}: Hopf algebras, topological vector spaces, 
quantum groups, deformation quantization.

\section{Introduction}
\subsection{Presentation of the context.}
The expression ``quantum groups'' is a name coined by Drinfeld (see 
\cite{Driv87a}) in the first half of the 80's which is superb, even if the 
notion is not necessarily quantum and the objects are not really groups. 
But they are Hopf algebras and their theory can be viewed as an avatar of 
deformation quantization \cite{BFFLS78} (see \cite{DS02} for a recent review 
which this presentation complements), applied to the quantization of 
Poisson-Lie groups. 

The philosophy underlying the role of {\it deformations in physics} 
has been consistently put forward by Flato, almost since the definition
of the deformation of rings and algebras by Gerstenhaber \cite{Ge64},
and was eventually expressed by him in \cite{Fl82}. 
In short, the passage from one level of physical theory to another,
when a new fundamental constant is imposed by experiments, can be 
understood (and might even have been predicted) using deformation theory. 
The only question is, in which category do we seek for deformations? 
Usually physics is rather conservative and if we start e.g. with the 
category of associative or Lie algebras, we tend to deform in the same 
category. 

But there are important instances of generalizations of this principle. 
The most elaborate is maybe noncommutative geometry, where the 
strategy is to formulate the ``undeformed" (commutative) geometry 
in terms of algebraic structures in such a way that it becomes possible 
to ``plug in" the deformation (noncommutativity) in a quite natural, 
and mathematically rigorous, manner. We shall not elaborate on that 
aspect here, refering e.g. to \cite{Co2k} for a presentation, 
to \cite{CDV02} for important recent examples of noncommutative manifolds,
and to \cite{Co94,CFS92} for the basics and a relation with 
deformation quantization. 

We shall concentrate on another prominent example: quantum groups.
Instead of looking at the associative algebra of functions over a 
Poisson-Lie group or at the enveloping algebra, one makes full use of 
the Hopf algebra structure in both cases. In general both the product and 
the coproduct have to be (compatibly) deformed, but
cohomological results (\cite{Driv89b} and section \ref{quantif})
show that, when the Lie group is semi-simple, the deformation is always
equivalent to a ``preferred" one, that is, 
a deformation where only the product or the coproduct (resp.) is deformed. 
The group aspect is a special case of deformation quantization and we 
shall show that the enveloping algebra aspect can be seen as its dual, 
in the sense of topological vector spaces duality.

\subsection{Deformation theory of algebras.}
A concise formulation of a Gerstenhaber deformation of an algebra 
(associative, Lie, bialgebra, etc.) is \cite{Ge64,BFGP94a}:
\begin{dfn}\label{defdef}
A deformation of an algebra $A$ over a field ${\mathbb K}$ is a 
${\mathbb K}[[\nu]]$-algebra ${\tilde A}$ such that 
${\tilde A}/\nu {\tilde A} \approx A$.
Two deformations ${\tilde A}$ and ${\tilde A'}$ are said equivalent if they
are isomorphic over ${\mathbb K}[[\nu]]$ and ${\tilde A}$ is said trivial if
it is isomorphic to the original algebra $A$ considered by base field
extension as a ${\mathbb K}[[\nu]]$-algebra.
\end{dfn}
Whenever we consider a topology on $A$, ${\tilde A}$ is supposed to be
topologically free. For associative (resp. Lie) algebras, 
Definition \ref{defdef} tells us that there exists a new product 
$\ast$ (resp. bracket $[\cdot,\cdot]$) 
such that the new (deformed) algebra is again associative (resp. Lie).
Denoting the original composition laws by ordinary product 
(resp. $\{\cdot,\cdot\}$) this means that, for $u,v\in A$
(we can extend this to $A[[\nu]]$ by ${\mathbb K}[[\nu]]$-linearity) we have:
\begin{eqnarray}
u\ast v &=& uv + \sum_{r=1}^\infty \nu^r C_r(u,v) \label{a}\\
\left[u, v\right] &=& \{u,v\} + \sum_{r=1}^\infty \nu^r B_r(u,v)\label{l}
\end{eqnarray}
where the $C_r$ are Hochschild 2-cochains and the $B_r$ (skew-symmetric)
Chevalley 2-cochains, such that for $u,v,w\in A$ we have
$(u\ast v) \ast w=u\ast (v\ast w)$ and ${\cal{S}}[[u,v],w]=0$, where
${\cal{S}}$ denotes summation over cyclic permutations. 

For a (topological) {\it bialgebra} (an associative algebra $A$ where we have
in addition a coproduct $\Delta : A \longrightarrow A \otimes A$ and the
obvious compatibility relations), denoting by $\otimes_\nu$ the tensor
product of ${\mathbb K}[[\nu]]$-modules, we can identify
${\tilde A}\, {\hat{\otimes}}_\nu {\tilde A}$ with
$(A\, {\hat{\otimes}}A)[[\nu]]$, where ${\hat{\otimes}}$ denotes the algebraic
tensor product completed with respect to some topology (e.g. projective 
for Fr\'echet nuclear topology on A), we similarly have a deformed
coproduct ${\tilde \Delta } = \Delta + \sum_{r=1}^\infty \nu^r D_r$,
$D_r \in {\cal L}(A, A {\hat{\otimes}}A)$, satisfying 
${\tilde \Delta }(u * v)={\tilde \Delta }(u) * {\tilde \Delta }(v)$.
In this context appropriate cohomologies can be introduced 
\cite{GeSh90a,Bonp92b}. 
There are natural additional requirements for Hopf algebras.
\medskip

{\it Equivalence} means that there is an isomorphism
$T_\nu=I+\sum_{r=1}^\infty \nu^r T_r$, $T_r\in{\cal L}(A,A)$ so that
$T_\nu(u*'v)=(T_\nu u*T_\nu v)$ in the associative case, denoting by
$*$ (resp. $*'$) the deformed laws in ${\tilde A}$ (resp. ${\tilde A'}$);
and similarly in the Lie, bialgebra and Hopf cases.
In particular we see (for $r=1$) that a deformation is trivial at order 1 if
it starts with a 2-cocycle which is a 2-coboundary. More generally, exactly
as above, we can show \cite{BFFLS78} (\cite{GeSh90a,Bonp92b} in 
the Hopf case) that if two deformations are equivalent
up to some order $t$, the condition to extend the equivalence one step
further is that a 2-cocycle (defined using the $T_k$, $k\leq t$) is the
coboundary of the required $T_{t+1}$ and therefore {\it the obstructions to
equivalence lie in the 2-cohomology}. In particular, if that space is null,
all deformations are trivial.
\medskip

{\it Unit}. An important property is that a {\it deformation of an
associative algebra with unit} (what is called a unital algebra) is
again unital, and {\it equivalent to a deformation with the same unit}.
This follows from a more general result of Gerstenhaber (for deformations
leaving unchanged a subalgebra) and a proof can be found in \cite{GS88}.
\medskip

\begin{rmk}\rm{ In the case of (topological) {\it bialgebras} or {\it Hopf}
algebras, {\it equivalence} of deformations has to be understood as an 
isomorphism of (topological) ${\mathbb K}[[\nu]]$-algebras, the isomorphism
starting with the identity for the degree 0 in $\nu$. A deformation is again
said {\it trivial} if it is equivalent to that obtained by base field
extension. For Hopf algebras the deformed algebras may be taken
(by equivalence) to have the same unit and counit, but in general not the
same antipode.}
\end{rmk}

\subsection{Deformation quantization and physics.}
Intuitively, classical mechanics is the limit of quantum mechanics 
when $\hbar=\frac{h}{2\pi}$ goes to zero. But how can this be realized 
when in classical mechanics the observables are functions over phase space 
(a Poisson manifold) and not operators? The deformation philosophy promoted 
by Flato shows the way: one has to look for deformations of algebras of 
classical observables, functions over Poisson manifolds, and realize there
quantum mechanics in an \textsl{autonomous} manner. 

What we call ``deformation quantization" relates to (and generalizes)
what in the conventional (operatorial) formulation are the Heisenberg 
picture and Weyl's quantization procedure. In the latter \cite{We31}, 
starting with a classical observable $u(p,q)$, some function on phase space
${\mathbb R}^{2\ell}$ (with $p, q \in {\mathbb R}^{\ell}$), one associates 
an operator (the corresponding quantum observable) $\Omega(u)$ in the
Hilbert space $L^2({\mathbb R}^{\ell})$ by the following general recipe:
\begin{equation} \label{weyl}
u \mapsto \Omega_w(u) = \int_{{\mathbb R}^{2\ell}} \tilde{u}(\xi,\eta)
{\exp}(i(P.\xi + Q.\eta)/\hbar)w(\xi,\eta)\ d^\ell \xi d^\ell \eta
\end{equation}
where $\tilde{u}$ is the inverse Fourier transform of $u$, $P_\alpha$ and
$Q_\alpha$ are operators satisfying the canonical commutation relations
$[P_\alpha , Q_\beta] = i\hbar\delta_{\alpha\beta}$
($\alpha, \beta = 1,...,\ell$), $w$ is a weight function and the integral
is taken in the weak operator topology. 
What is called in physics normal (or antinormal) ordering corresponds to 
choosing for weight $w(\xi,\eta) = {\exp}(-{1\over 4}(\xi^2 \pm \eta^2))$.
Standard ordering (the case of the usual pseudodifferential operators in
mathematics) corresponds to $w(\xi,\eta) = {\exp}(-{i\over 2}\xi\eta)$ and 
the original Weyl (symmetric) ordering to $w = 1$.
An inverse formula was found shortly afterwards by Eugene Wigner \cite{Wi32}
and maps an operator into what mathematicians call its symbol by a kind
of trace formula. For example $\Omega_1$ defines an isomorphism of Hilbert
spaces between $L^2({\mathbb R}^{2\ell})$ and Hilbert-Schmidt operators on
$L^2({\mathbb R}^{\ell})$ with inverse given by
\begin{equation} \label{EPW}
u=(2\pi\hbar)^{-\ell}\, {\rm{Tr}}[\Omega_1(u)\exp((\xi.P+\eta.Q)/i\hbar)]
\end{equation}
and if $\Omega_1(u)$ is of trace class one has
${\rm{Tr}}(\Omega_1(u))=(2\pi\hbar)^{-\ell}\int u \, \omega^\ell
\equiv \mathrm{Tr}_\mathrm{\scriptstyle M}(u)$, the ``Moyal trace", where
$\omega^\ell$ is the (symplectic) volume $dx$ on ${\mathbb R}^{2\ell}$.
Looking for a direct expression for the symbol of a quantum commutator,
Moyal found \cite{Mo49} what is now called the Moyal bracket:
\begin{equation} \label{Moyal}
M(u_1,u_2) = \nu^{-1} \sinh(\nu P)(u_1,u_2) = P(u_1,u_2) + 
\sum^\infty_{r=1}\frac{\nu^{2r}}{(2r+1)!} P^{2r+1} (u_1,u_2) 
\end{equation}
where $2\nu=i\hbar$, $P^r(u_1,u_2)=\Lambda^{i_1j_1}\ldots
\Lambda^{i_rj_r}(\partial_{i_1\ldots i_r}u_1)(\partial_{j_1\ldots j_r} u_2)$
is the $r^{\rm{th}}$ power ($r\geq1$) of the Poisson bracket
bidifferential operator $P$, $i_k, j_k = 1,\ldots,2\ell$, $k=1,\ldots,r$
and $(\Lambda^{i_kj_k}) = {0\,-I\choose I\,0}$. To fix ideas we may 
assume here $u_1,u_2\in \mathcal{C}^\infty({\mathbb R}^{2\ell})$ and 
the sum is taken as a formal series. 
A corresponding formula for the symbol of a product $\Omega_1(u)\Omega_1(v)$ 
can be found in \cite{Gr46}, and may now be written more clearly as 
a (Moyal) {\it star product}:
\begin{equation} \label{star}
u_1 \ast_M u_2 = \exp(\nu P)(u_1,u_2) = 
u_1u_2 + \sum^\infty_{r=1}\frac{\nu^{r}}{r!} P^{r}(u_1,u_2).
\end{equation}
The formal series may be deduced (see e.g. \cite{Biep00a}) from an integral 
formula of the type:
\begin{equation}\label{starint}
(u_1\ast u_2)(x)=c_\hbar \int_{\mathbb{R}^{2\ell}\times\mathbb{R}^{2\ell}}
u_1(x+y)u_2(x+z)e^{-\frac{i}{\hbar}\Lambda^{-1}(y,z)}dydz.
\end{equation}
It was noticed, however after deformation quantization was introduced, 
that the composition of symbols of pseudodifferential operators (ordered, 
like differential operators, ``first $q$, then $p$") is a star product. 

One recognizes in (\ref{star}) a special case of (\ref{a}), and 
similarly for the bracket. So, via a Weyl quantization map, the algebra 
of quantized observables can be viewed as a deformation of that of classical
observables. 

But the deformation philosophy tells us more.
Deformation quantization is not merely ``a reformulation of quantizing 
a mechanical system'' \cite{DoNe01}, e.g. in the framework of Weyl 
quantization: \textsl{The process of quantization itself is a deformation}.
In order to show that explicitly it was necessary to treat in an 
\textsl{autonomous} manner significant physical examples, without recourse 
to the traditional operatorial formulation of quantum mechanics. 
That was achieved in \cite{BFFLS78} with the paradigm of the harmonic 
oscillator and more, including the angular momentum and the hydrogen atom. 
In particular what plays here the role of the unitary time evolution operator 
of a quantized system is the ``star exponential" of its classical Hamiltonian 
$H$ (expressed as a usual exponential series but with ``star powers" of 
$tH/i\hbar$, $t$ being the time, and computed as a distribution both in phase
space variables and in time); in a very natural manner, the spectrum 
of the quantum operator corresponding to $H$ is the support of the 
Fourier-Stieltjes transform (in $t$) of the star exponential (what Laurent
Schwartz had called the spectrum of that distribution). Further examples 
were (and are still being) developed, in particular in the direction of
field theory.

That aspect of deformation theory has since 25 years or so been extended 
considerably. It now includes general symplectic and Poisson 
(finite dimensional) manifolds, with further results for infinite dimensional 
manifolds, for ``manifolds with singularities" and for algebraic varieties, 
and has many far reaching ramifications in both mathematics and physics 
(see e.g. a brief overview in \cite{DS02}). As in quantization 
itself \cite{We31}, symmetries (group theory) play a special role 
and an autonomous theory of star representations of Lie groups 
was developed, in the nilpotent and solvable cases of course (due to the 
importance of the orbit method there), but also in significant other examples.
The presentation that follows can be seen as an extension of the latter,
when one makes full use of the Hopf algebra structures and of the ``duality"
between the group structure and the set of its irreducible representations.

Finally one should mention that deformation theory and Hopf algebras
are seminal in a variety of problems ranging from theoretical physics
(see e.g. \cite{CK99,DS02}), including renormalization and Feynman integrals 
and diagrams, to algebraic geometry and number theory (see e.g.
\cite{Ko01,KZ01}), including algebraic curves \`a la Zagier 
(cf. \cite{CM03} and Connes' lectures at Coll\`ege de France, 
January to March 2003). 

\section{Some topological Hopf algebras}

We shall now briefly review applications of the deformation theory of 
algebras in the context of Hopf algebras endowed with appropriate 
topologies and in the spirit of deformation quantization. 
That is, we shall consider Hopf algebras of functions on Poisson-Lie groups 
(or their topological duals) and their deformations, 
and show how this framework is a powerful tool to understand the 
standard examples of quantum groups, and more. 
In order to do so we first recall some notions on topological
vector spaces and apply them to our context.

\subsection{Well-behaved Hopf algebras}   
\begin{dfn}
A topological vector space (tvs) $V$ is said {\em well-behaved} 
if $V$ is either nuclear and Fr\'echet, or nuclear and dual of Fr\'echet
\cite{Groa55, Tref67}.
\end{dfn}

\begin{prop}
If $V$ is a well-behaved tvs and $W$ a tvs, then 
$$(i)\ V^{**}\simeq V \qquad 
(ii)\ (V\hat{\otimes}V)^*\simeq V^*\hat{\otimes}V^* \qquad
(iii)\ \Hom_\mathbb{K}(V,W)\simeq V^*\hat{\otimes} W$$
where $V^*$ denotes the strong topological dual of $V$, $\hat{\otimes}$
the projective topological tensor product and the base field
$\mathbb{K}$ is $\mathbb{R}$ or $\mathbb{C}$. 
\end{prop}

\begin{dfn}
$(A,\mu,\eta,\Delta,\epsilon,S)$ is a WB (well-behaved) Hopf 
algebra \cite{BFGP94a} if
\begin{itemize}
\item $A$ is a well-behaved topological vector space.
\item The multiplication $\mu :A\hat{\otimes}A\rightarrow A\ $, 
the coproduct $\Delta : A\rightarrow A\hat{\otimes}A\ $, the unit $\eta\ $, 
the counit $\epsilon\ $, and the antipode $S$ are continuous. 
\item $\mu , \eta , \Delta , \epsilon$ and $S$ satisfy the usual axioms
of a Hopf algebra.
\end{itemize}
\end{dfn}

\begin{cor}
If $(A, \mu , \eta , \Delta , \epsilon , S)$ is a WB Hopf algebra, then
$(A^*,\,^t\Delta ,\,^t\epsilon ,\,^t\mu ,\,^t\eta ,\,^tS)$ is also 
a WB Hopf algebra.
\end{cor}

\subsection{Examples of well-behaved Hopf algebras \cite{BFGP94a}}

Let $G$ be a semi-simple Lie group and $\mathfrak{g}$ its complexified
Lie algebra. For simplicity we shall assume here $G$ linear (i.e. with a 
faithful finite dimensional representation) but the same results hold, with 
some modification in the proofs, for any semi-simple Lie group.

\subsubsection{Example 1}

$\mathcal{C}^\infty (G)$, the algebra of the smooth functions on $G$,
is a WB Hopf algebra (Fr\'echet and nuclear).

\subsubsection{Example 2}

$\mathcal{D}(G)=\mathcal{C}^\infty (G)^*$, the algebra of the compactly 
supported distributions on $G$, is a WB Hopf algebra (dual of Fr\'echet 
and nuclear). The product is the transposed map of the coproduct of 
$\mathcal{C}^\infty (G)$ that is, the convolution of distributions.

\subsubsection{Example 3}

$\mathcal{H}(G)$, the algebra of coefficient functions of finite
dimensional representations of $G$ (or polynomial functions on $G$)
is a WB Hopf algebra, the Hopf structure being that induced from 
$\mathcal{C}^\infty (G).$

A short description of that algebra is as follows: We take a set $\hat{G}$ 
of irreducible finite dimensional representations of $G$ such that
there is \textsl{one and only one} element for each equivalence class,
and, if $\pi \in \hat{G}$, its contragredient $\check{\pi}$ is also in
$\hat{G}$. We define
$\ds \quad C_\pi = \mathsf{vect}\{\mbox{coefficient functions of } \pi\}
\stackrel{Burnside}{\simeq} \End (V_\pi )$ for $\pi \in \hat{G}$. 
Then $\ds \quad \mathcal{H}(G)
\stackrel{alg.}{\simeq} \bigoplus_{\pi \in \hat{G}}C_\pi 
\stackrel{v. s.}{\simeq} \bigoplus_{\pi \in \hat{G}}\End (V_\pi ).$
So we take on $\mathcal{H}(G)$ the ``direct sum'' topology of 
$\ds \bigoplus_{\pi \in \hat{G}}\End (V_\pi )$. Then $\mathcal{H}(G)$
is dual of Fr\'echet and nuclear, that is, WB.

\subsubsection{Example 4}

Let $\mathcal{A}(G)$, the algebra of ``generalized distributions'', 
be defined by $ \mathcal{A}(G)= \mathcal{H}(G)^*
\stackrel{alg.}{\simeq}\prod_{\pi \in \hat{G}}\End (V_\pi ).$ 
The (product) topology is Fr\'echet and nuclear, and therefore 
$\mathcal{A}(G)$ is WB.

\subsection{Inclusions \cite{BiPi96, BFGP94a}}

We denote by $\mathsf{U}\mathfrak{g}$ the universal enveloping algebra
of $\mathfrak{g}$ and by $\mathbb{C}G$ the group algebra of $G.$
All the following inclusions are inclusions of Hopf algebras.
$\Subset$, $\Supset$, $\Cup$, $\Cap$ mean a {\em dense} inclusion.

$$
\begin{array}{cccccc|ccc}
\mathsf{U}\mathfrak{g} & \Subset & \mathcal{A}(G) 
& \Supset & \mathbb{C}G & & & \mathcal{H}(G) & \\
& & \Cup & & & & & \Cap & (*) \\
\mathsf{U}\mathfrak{g} & \subset & \mathcal{D}(G) 
& \Supset & \mathbb{C}G & & & \mathcal{C}^\infty (G) &
\end{array}$$

\hfill $(*)$ is true if and only if $G$ is linear, 
but comparable results can be obtained for $G$ non linear.

\section{Topological quantum groups}
We shall now deform the preceding topological Hopf algebras and indicate
how this explains various models of quantum groups. For clarity of the 
exposition, throughout this Section and the remainder of the paper, we
shall limit to a minimum the details concerning the Hopf algebra structures
other than product and coproduct. But whenever we write Hopf algebras 
and not only bialgebras, the relevant structures are included in the 
discussion and dealing with them is quite straightforward.

\subsection{Quantization }\label{quantif}

\begin{thm}[\cite{Driv89b}] \label{rigid}
Let $\mathfrak{g}$ be a semi-simple Lie algebra and 
$ (\mathsf{U}\mathfrak{g} , \mu_0 , \Delta_0 ) $ denote
the usual Hopf structure on $\mathsf{U}\mathfrak{g}$.
\begin{enumerate}
\item If $\ds\ (\mathsf{U}_t\mathfrak{g} , \mu_t ) $ is a deformation
(as an algebra) of $\ds\ (\mathsf{U}\mathfrak{g}[[ t ]] , \mu_0 )$
then $\ds\ \mathsf{U}_t\mathfrak{g}\stackrel{\varphi}{\simeq} 
\mathsf{U}\mathfrak{g}[[ t ]]\ $ 
(i.e. $\mathsf{U}\mathfrak{g}$ is rigid). 
\item If $\ds\ (\mathsf{U}\mathfrak{g}[[ t ]] , \mu_0 , \Delta_t ) $
is a deformation (as a Hopf algebra) of
$\ds\ (\mathsf{U}\mathfrak{g}[[ t ]] , \mu_0 , \Delta_0 ) $ 
then $$\exists \ P_t\in (\mathsf{U}\mathfrak{g}\otimes\mathsf{U}\mathfrak{g})
[[ t ]] \ such\ that\ P_{t=0}=\Id \ \mbox{ and }\ 
\Delta_t(a)=P_t . \Delta_0(a) . P_t^{-1},
\ \forall a\in \mathsf{U}\mathfrak{g}.$$
\end{enumerate}
\end{thm}

An isomorphism $\varphi$ (it is not unique!) appearing in item 1 above 
is called a {\em Drinfeld isomorphism}.

\begin{cor}[\cite{BFGP94a}]
Let $G$ be a linear semi-simple Lie group and $\mathfrak{g}$ be its 
complexified Lie algebra. 
\begin{enumerate}
\item If $\mathsf{U}_t\mathfrak{g}$ is a deformation of 
$\mathsf{U}\mathfrak{g}$ (a ``quantum group'') then
$\ds (\mathsf{U}_t\mathfrak{g},\mu_t,\Delta_t)
\simeq (\mathsf{U}\mathfrak{g} [[ t ]],\mu_0, P_t \Delta_0 P^{-1}_t).$
\item $\ds\mathcal{ A}_t(G):=(\mathcal{ A} (G)[[ t ]],
\mu_0,P_t\cdot\Delta_0\cdot P_t^{-1})$ 
is a Hopf deformation of $\mathcal{ A} (G)$
and $\ds \mathsf{U}_t\mathfrak{g}
\stackrel{\mathrm{Hopf}}{\subset}\mathcal{ A}_t (G).$ 
\item $\ds \mathcal{ D}_t(G):=(\mathcal{ D}(G)[[ t ]],
\mu_0, P_t\cdot\Delta_0\cdot P^{-1}_t)$ 
is a Hopf deformation of $\mathcal{ D}(G)$ and 
$\ds\mathsf{U}_t\mathfrak{g}\stackrel{\mathrm{Hopf}}{\subset}\mathcal{D}_t(G)$.
\item $\ds\mathcal{ C}^\infty_t(G):=\mathcal{ D}_t(G)^\ast$ and 
$\ds\mathcal{ H}_t(G):=\mathcal{ A}_t(G)^\ast$ are
quantized algebras of functions. They are Hopf deformations of 
$\mathcal{ C}^\infty(G)$ and $\mathcal{ H}(G).$
\end{enumerate}
\end{cor}
Similar results hold in the non linear case \cite{BiPi96} and for other WB
Hopf algebras (e.g. constructed with infinite dimensional representations)
\cite{Bidf96a}.

\Pf
(1) Direct consequence of Theorem \ref{rigid}.\ \ 
(2) $\ds P_t\in (\mathsf{U}\mathfrak{g}\otimes 
\mathsf{U}\mathfrak{g})[[ t ]]\subset 
(\mathcal{ A} (G)\hat\otimes\mathcal{ A} (G))[[ t ]]$.
We obtain coassociativity from 
$\mathsf{U}\mathfrak{g} \Subset \mathcal{ A} (G).$
(3) By restriction of (2).
(4) By simple dualization from (2) and (3).
\EPf

\begin{rmk} ``Hidden group structure" in a quantum group. 
{\rm Here the deformations are {\em preferred}, that is, the product
on $\mathcal{ D}_t(G)$ and on $\mathcal{ A}_t(G)$ (resp. the coproduct
on $\mathcal{ C}^\infty_t(G)$ and on $\mathcal{ H}_t(G)$) is not deformed
and the basic structure is still the product on the group $G$.
So this approach gives an interpretation of the Tannaka-Krein philosophy 
in the case of quantum groups: it has often been noticed that, in the
generic case, finite dimensional representations of a quantum group are
(essentially) representations of its classical limit. So the algebras 
involved should be the same, which is justified by the above mentioned 
rigidity result of Drinfeld. This shows that the initial classical group 
is still there, acting as a kind of ``hidden variables" in this quantum 
group theory, which is exactly what we see in this quantum group theory.
This fact was implicit in Drinfeld's work. The Tannaka-Krein interpretation
of the twisting of quasi-Hopf algebras can be found in Majid
(see e.g. \cite{Majs92a}). It was made explicit, within the framework
exposed here, in \cite{BFGP94a}.}
\end{rmk}

\subsection{Unification of models and generalizations}

\subsubsection{Drinfeld models}
We call ``Drinfeld model of quantum group" a deformation of
$\mathsf{U}\mathfrak{g}$ for $\mathfrak{g}$ simple, as given in
\cite{Driv87a}. We have seen in the preceding section that from 
any Drinfeld model $\mathsf{U}_t\mathfrak{g}$ of a quantum group
(which can be generalized to any deformation of the Hopf algebra 
$\mathsf{U}\mathfrak{g}$), we obtain a deformation of 
$\mathcal{D}(G)$ and $\mathcal{A}(G)$ that contains 
$\mathsf{U}_t\mathfrak{g}$ as a sub-Hopf algebra. 
So $\mathcal{D}_t(G)$ and $\mathcal{A}_t(G)$ are quantum group models 
that describe Drinfeld models. By duality, $\mathcal{C}_t^\infty (G)$ 
and $\mathcal{ H}_t(G)$ are ``quantum group deformations" of
$\mathcal{C}^\infty (G)$ and $\mathcal{ H} (G)$. 
The deformed product on $\mathcal{H}(G)$ is the restriction of that
on $\mathcal{C}^\infty (G)$. Furthermore, as we shall see, 
these deformations coincide with the usual ``quantum algebras of functions''. 
Let us look more in detail at $\mathcal{ H}_t(G)$:  

\subsubsection{Faddeev-Reshetikhin-Takhtajan (FRT) models}

In \cite{FRT88a} quantized algebras of functions are defined in terms 
of generators and relations, the key relation being given by the 
star-triangle (Yang-Baxter) equation, 
$\ R(T\otimes \Id)(\Id \otimes T)=(\Id \otimes T)(T\otimes \Id)R\ ,$
for a given R-matrix $R\in \End (V\otimes V)$ and for $T\in\End (V),$
$V$being a finite dimensional vector space.

\smallskip

As our deformations are given by a twist $P_t$, it is not surprising,
from a structural point of view \cite{Majs92a} that, dually, we
obtain in each case a Yang-Baxter relation and so a ``FRT-type''
quantized algebra of functions. Our Fr\'echet-topological context 
permits to write precisely such a construction for the infinite-dimensional
Hopf algebras involved.

\medskip

\textsl{3.2.2.1. Linear case.} If $G$ is semi-simple and linear, there exists 
$\pi$ a finite dimensional representation of $G$ such that 
$\mathcal{H}(G)\simeq \mathbb{ C}[\pi_{ij};1\leqslant i,j\leqslant N]$ 
where the $\pi_{ij}$ are the coefficient functions of $\pi$.
Denote by $(\mathcal{ H}_t (G),\ast)$ the deformation of $\mathcal{ H} (G)$ 
obtained in this way and by $T$ the matrix $[\pi_{ij}]$. 
Define $T_1:=T\otimes Id$ and $T_2:= Id\otimes T$. Then we have

\begin{prop}[\cite{BFGP94a, BiPi96}] \hfill

\begin{enumerate}
\item $\{ \pi_{ij};1\leqslant i,j\leqslant N\}$is a topological 
generator system of the $\mathbb{C}[[t]]$-algebra $\mathcal{H}(G)_t$.
\item There exists an invertible $\mathcal{ R}\in \mathcal{ L}
(V_\pi\otimes V_\pi)[[ t ]]$ such that
$\mathcal{ R}\cdot T_1\ast T_2=T_2\ast T_1\cdot \mathcal{ R}$
(so $\mathcal{ H}_t (G)$ is a ``quantum algebra of functions" of type FRT).
\item We recover every quantum group 
given in \cite{FRT88a} by this construction.
\end{enumerate}
\end{prop}
\SPf
\begin{enumerate}
\item Perform a precise study of the deformed tensor product of representations.
\item  Since the deformations $\mathcal{A}_t(G)$ are given by a 
twist $P_t$, $\mathcal{A}_t(G)$ is quasi-cocommutative, i.e. there exists 
$ R \in (\mathcal{A}(G)\hat{\otimes}\mathcal{A}(G))[[ t ]]$
such that $\sigma \circ \Delta_t (a) = R \Delta_t(a) R^{-1}$ with
$\sigma (a\otimes b)= b\otimes a$.
Standard computations give the result.

\item We want to follow the way used in \cite{Driv87a} to link Drinfeld 
to FRT models. But the main point is that our deformations are obtained 
through a Drinfeld isomorphism. We therefore have to show:

- There exists a specific Drinfeld isomorphism deforming the
standard representation of $\mathfrak{g}$ into the representation
of $\mathsf{U}_t\mathfrak{g}$ used in \cite{Driv87a}.

- Two Drinfeld isomorphisms give equivalent deformations. \hfill \EPf
\end{enumerate}

For instance, the FRT quantization of $SL(n)$ can be seen as a Hopf 
deformation of $\mathcal{ H}(SU(n))$ (with non deformed coproduct).
Moreover, this Hopf deformation extends to $\mathcal{C}^\infty (G)$.
\newpage

\begin{rmk}\hfill
\begin{enumerate}
\item {\rm This proposition justifies the terminology ``deformation'', 
often employed but never justified in these cases. 
See e.g. \cite{GGS91} where it is shown that relations of type
$\mathcal{ R} T_1 T_2=T_2 T_1 \mathcal{ R}$ need
not define a deformation, even if $\mathcal{ R}$ is Yang-Baxter.}
\item {\rm Starting from Drinfeld models, our construction produces FRT models 
also for e.g. $G=Spin(n)$ and for exceptional Lie groups.
In addition, at least some multiparameter deformations \cite{Res90}
can be easily treated in this way \cite{BFGP94a}.}
\end{enumerate}
\end{rmk} 

\noindent \textsl{3.2.2.2. Non-linear case.}

\begin{prop}[\cite{BiPi96}]
If $G$ is semi-simple with finite center, there exists 
a dense subalgebra of $(\mathcal{C}^\infty_t(G),\ast)$
generated by the coefficient functions of a finite number 
of (possibly infinite dimensional) representations.
\end{prop}

\subsubsection{Jimbo models}  

These are models \cite{Jim85} with generators $E_i^\pm$, $K_i$ and $K^{-1}_i$.
For $G=SU(2)$ \cite{BFP92a} and $G=SL(2,\mathbb{C})$ \cite{MaZo96} 
we realize $\mathsf{U}_q \mathfrak{sl}(2)$ and 
$\mathsf{U}_t\mathfrak{sl}(2,\mathbb{C})$ as dense sub-Hopf algebras of 
$\mathcal{A} (G)$, $\forall t\in \mathbb{C}\setminus 2\pi \mathbb{ Q}$
(with $q=e^t$). For $\mathfrak{sl}(2)$ this gives the original model of 
Jimbo \cite{Jim85}. For the Lorentz algebra $\mathfrak{sl}(2,\mathbb{C})$ 
this unifies \cite{MaZo96} all the models proposed so far in the literature 
for a quantum Lorentz group. We obtain here {\em convergent} deformations 
(not only formal). 

For $\mathfrak{sl}(2,\mathbb{C})$ it was first proposed in \cite{PW90} 
to consider the quantum double \cite{Driv87a} of 
$\mathsf{U}_q \mathfrak{su}(2)$ as $q$-deformed Lorentz group. It was known
from \cite{RSts90} that in such cases the double, as an algebra, is
the tensor product of two copies of $\mathsf{U}_t\mathfrak{su}(2)$.
See also \cite{OSWZ91a,SWZ91a}, and \cite{Maj93} for a dual version and
another semi-direct product form. 

\subsubsection{Deformation quantization}

{}From the main construction, using deformations of $\mathsf{U}\mathfrak{g},$ 
we deduce the following general theorem:

\begin{thm}[\cite{BiPi96}] Let $G$ be a semi-simple connected Lie group 
with a Poisson-Lie structure. There exists a deformation
$(\mathcal{C}^\infty_t (G),\ast)$ of $\mathcal{C}^\infty (G)$ such that
$\ast$ is a (differential) star product.
\end{thm}

\begin{rmk}\hfill
{\rm \begin{itemize}
\item When Lie$(G)$ is the double of some Lie algebra, the same result holds.
\item The fact that $\ast$ is differential comes from the twist 
$P_t\Delta_0P^{-1}_t$, 
$P_t\in (\mathsf{U}\mathfrak{g}\times \mathsf{U}\mathfrak{g})[[ t ]]$.
\item{Since from any Drinfeld quantum group we obtain a star product, and 
since any FRT quantum group can be seen as a restriction of such a 
star product, we have showed that the data of a ``semi-simple'' quantum group 
is equivalent to the data of a star product on $\mathcal{C}^\infty(G)$ 
satisfying $\Delta (f\ast g)=\Delta (f)\ast\Delta (g)$.
The functorial existence results of Etingof and Kazhdan \cite{EtKa96a} on the 
quantization of Lie bialgebras (see also \cite{Enr02}) show that
the latter is true also for ``non semi-simple" quantum groups.}
\item{Techniques similar to those indicated here can be applied to other 
$q$-algebras (more general quantum groups such as those in \cite{Fro97}
and more recent examples, Yangians, etc.). 
In particular those used in the case of the Jimbo models should be applicable
to $q$-algebras defined by generators and relations. That direction of
research has not yet been developed.}
\end{itemize}}
\end{rmk}

\section{Topological quantum double}

{}From now on we use the Sweedler notation for the coproducts \cite{Swem68a}:
in a coalgebra $(H,\Delta)$, $\Delta(x)=\sum_{(x)}x_{(1)}\otimes x_{(2)}$
and, by coassociativity, 
$(\Id\otimes\Delta)\Delta(x) = (\Delta\otimes\Id)\Delta(x) 
= \sum_{(x)}x_{(1)}\otimes x_{(2)}\otimes x_{(3)}.$

\medskip

In \cite{Driv87a} Drinfeld defines the quantum double of 
$\mathsf{U}_t\mathfrak{g}$ (see also \cite{Semm94a}). This can be adapted 
to the context of topological Hopf algebras \cite{Bonp94b}.

\subsection{Definitions}

Let $A$ be $\mathcal{ D} (G),\mathcal{ A} (G),\mathcal{ D}_t(G)$ 
or $\mathcal{ A}_t(G).$ If $A=(A,\mu,\Delta, S)$ then
$\ds A^\ast = (A^\ast,\,^t\Delta,\,^t\mu,\,^tS)$. Define $A^0=A^{\ast~co-op}
=(A^\ast,\,^t\Delta,\,^t\mu^{op},\,^tS^{op}),$
where $\mu^{op}(x\otimes y):= \mu (y\otimes x)$ and $S^{op}$ is 
the antipode compatible with $\mu^{op}$ and $\Delta .$

\smallskip

If we consider the vector space $A^* \otimes A,$ Drinfeld \cite{Driv87a} 
defines the quantum double as follows :

i) $D(A) \simeq A^0 \otimes A$ as coalgebras,

ii) $(f \otimes Id_A) . (Id_{A^0} \otimes b) = f \otimes b,$

iii) $(Id_{A^0} \otimes e_s) . (e^t \otimes Id_A) = \Delta^{kjn}_s \hskip.1cm
\mu^t_{plk} \hskip.1cm S'^p_n \hskip.1cm (e^l \otimes Id_A) \hskip.1cm 
(Id_{A^0} \otimes e_j),$ where $\{ e_s \}$ is a basis of $A$ and 
$\hskip.1cm \{ e^t \}$ the dual basis.

The Drinfeld double was expressed \cite{Majs90a} in a Sweedler form for 
dually paired Hopf algebras as an example of a theory of `double smash 
products'. Adapting that formulation to our topological context we can now 
define the double as:

\begin{dfn}
The double of $A$, $D(A)$, is the topological Hopf algebra 
$\ds (A^\ast\overline\otimes A,\mu_D,\,^t\mu^{op}
\otimes\Delta,\,^tS^{op}\otimes S)$ with 
\begin{eqnarray*} \ds \mu_D((f\otimes a)\otimes (g\otimes b)) & = &
\sum_{(a)} f<g\,,\, S^{op}\,(a_{(3)})\,?\,a_{(1)}> \otimes a_{(2)}\,b\\
 & = & \sum_{(a) (g)} <g_{(1)}, a_{(1)}>
\hskip.1cm <\,^tS^{op}(g_{(3)}), a_{(3)}> \hskip.1cm fg_{(2)}\otimes a_{(2)} b
\end{eqnarray*} 
where $<~~,~~>$ denotes the pairing $A^\ast/ A,$
$``?"$ stands for a variable in $A$ and $\overline\otimes$ is the completed 
inductive tensor product.
\end{dfn}

As topological vector spaces we have $D(A)=A^\ast\overline\otimes A$. Thus 
$D(A)^\ast=A\hat\otimes A^\ast$ and $D(A)^{\ast\ast}=D(A)$.
So $D(A)$ is ``almost self dual'' (it is self dual up to a completion) 
and is reflexive.

\subsection{Extension theory}

\begin{itemize}
\item If $A$ is cocommutative then the product $\mu_D$ of 
$D(A)$ is the {\em smash product} $\buildrel\mu_{}^{\rightharpoonup}$ on 
$A^0\overline\otimes A$
$$\buildrel\mu_{}^{\rightharpoonup}\bigl((f\otimes a)\otimes (g\otimes b)\bigl)
=\sum_{(a)} f (a_{(1)}\rightharpoonup g)\otimes a_{(2)}b$$
where $\rightharpoonup$ denotes the coadjoint action of $A$ on $A^0$,
$<a\rightharpoonup f , b> = \sum_{(a)} <f , S(a_{(1)}) b a_{(2)} >$.
This product is the ``zero class" of an extension theory, defined 
by Sweedler \cite{Swem68a}, classified by a space of 2-cohomology 
$H^2_{sw}(A,A^0)$. The products are of the form, for $\tau$ a 2-cocycle,
$$\buildrel{\mu_{\tau}}_{}^{\rightharpoonup} \bigl((f\otimes a)
\otimes (g\otimes b)\bigl)=\sum_{(a)(b)} f(a_{(1)}\rightharpoonup
g)\tau (a_{(2)}\otimes b_{(2)})\bigl)\otimes a_{(3)}b_{(2)}.$$
\item The coproduct of $D(A)$ is a smash coproduct 
for the trivial co-action. We can dualize the theory and, putting the two
things together, we obtain an extension theory for bialgebras which is 
classified by a cohomology space $H^2_{bisw} (A^0,A)$.
\end{itemize}

\noindent \textsl{Question} : Are there other possible definitions 
of the double as an extension of $A^0$ by $A$?\\
\noindent \textsl{Answer} : NO, for $A=\mathcal{ D} (G)$ \cite{Bonp94b},
because $H^2_{bisw}\bigl(\mathcal{ D}(G),\mathcal{C}^\infty (G)\bigl)=\{0\}.$

\section{Crossed products and deformation quantization}

In this section we shall see that the Hopf algebra techniques presented 
in the preceding sections can be useful not only to understand quantum
groups, but also to develop very nice formulas in deformation quantization
itself.

In order to shed light on the general definition which follows, 
we return to the simplest case of deformation quantization: 
the Moyal product on $\mathbb{R}^2$. We look at $\mathbb{R}^2$ as 
$T^*\mathbb{R}\equiv \mathbb{R}\times\mathbb{R}^*$ and therefore can write
$\mathcal{C}^\infty(\mathbb{R}^2)\simeq \mathcal{C}^\infty(\mathbb{R})
\hat\otimes\mathcal{C}^\infty(\mathbb{R}^*).$
We consider first two functions of a special kind in this algebra:
$u(x)=u(x_1,x_2)=f(x_1)P(x_2)$ and $v(x)=v(x_1,x_2)=g(x_1)Q(x_2)$
where $f, g\in \mathcal{C}_0^\infty(\mathbb{R})$ and 
$P, Q$ are polynomials in $\Pol(\mathbb{R}^*)\simeq \mathsf{S}\mathbb{R}$.
We can then write is the usual coproduct on the symmetric
algebra $\mathsf{S}\mathbb{R}$ as 
$\ds\Delta(P)(x_2,y_2)=P(x_2+y_2) 
(\stackrel{\mathrm{\scriptscriptstyle notation}}{=} 
\sum_{(P)}P_{(1)}(x_2)P_{(2)}(y_2) )$.

We now look at Formula (\ref{starint}) for the Moyal star product on 
$\mathbb{R}^2$ and perform on it some formal calculations (we do not 
discuss the convergence of the integrals involved).
Up to a constant (depending on $\hbar$) we get:
\begin{eqnarray*}
(u\ast v)(x) & = &\int_{\mathbb{R}^{2}\times\mathbb{R}^{2}}
u(x+y)v(x+z)e^{-\frac{i}{\hbar}\Lambda^{-1}(y,z)}dydz\\
& = & \int_{\mathbb{R}^{2}\times\mathbb{R}^{2}}
f(x_1+y_1)P(x_2+y_2)g(x_1+z_1)Q(x_2+z_2)
e^{-\frac{i}{\hbar}(y_1z_2 - y_2z_1)}dy_1dy_2dz_1dz_2\\
& = & \int_{\mathbb{R}^{2}}f(x_1+y_1)Q(x_2+z_2)
e^{-\frac{i}{\hbar}y_1z_2}dy_1dz_2 . \int_{\mathbb{R}^{2}}
g(x_1+z_1)P(x_2+y_2)e^{\frac{i}{\hbar}y_2z_1}dy_2dz_1\\
& = & \sum_{(P)(Q)}(\partial^+_{Q_{(1)}}f)(x_1)Q_{(2)}(x_2) . 
(\partial^-_{P_{(1)}}g)(x_1)P_{(2)}(x_2)
\qquad \mbox{ (up to a constant) }
\end{eqnarray*}
with $\partial^\pm_{Q_{(1)}}=Q_{(1)}(\mp i\hbar \partial_{x_1})$ 
(the same for $P$), since 
$F_\hbar^\mp\left(\alpha F_\hbar^\pm (h)(\alpha)\right)(x)
= \mp i\hbar \partial_xh(x)$ for $h\in\mathcal{C}_0^\infty(\mathbb{R})$ with
$F_\hbar^\pm(h)(\alpha)$ defined as 
$\int_\mathbb{R}h(x)e^{\mp\frac{i}{\hbar}x\alpha}dx$.
This suggests the following small generalization of the smash product:

\begin{dfn}
Let $B$ be a cocommutative bialgebra and $C$ a $B$-bimodule algebra 
[i.e. $C$ is both a left $B$-module algebra and a right $B$-module algebra 
such that 
$(a\rightharpoonup f)\leftharpoonup b = a\rightharpoonup(f\leftharpoonup b)$].
We define the L-R smash product on $C\otimes B$ by 
$$(f\otimes a)\star (g\otimes b)=\sum_{(a)} 
(f\leftharpoonup b_{(1)})(a_{(1)}\rightharpoonup g)\otimes a_{(2)}b_{(2)}.$$
\end{dfn}

\begin{prop}
The L-R smash product is associative.
\end{prop}

\subsection{Relation with usual deformation quantization}

Let $G$ be a Lie group, $T^\ast G$ its cotangent bundle, 
$\mathfrak{g}=\hbox{Lie}(G)$. We have 
$$\mathcal{C}^\infty (T^\ast G)
\simeq \mathcal{C}^\infty (G\times \mathfrak{g}^\ast)
\simeq \mathcal{C}^\infty (G)\hat\otimes \mathcal{C}^\infty (\mathfrak{g}^\ast)
\supset\mathcal{C}^\infty (G)\otimes \mathsf{Pol}(\mathfrak{g}^\ast)
\simeq \mathcal{C}^\infty (G)\otimes \mathsf{S}\mathfrak{g}.$$
We define a deformation of 
$\mathcal{C}^\infty (G)\otimes \mathsf{S}\mathfrak{g}$ 
by a L-R smash product:

\begin{itemize}
\item We deform $\mathsf{S}\mathfrak{g}$ by the ``parametrized version" 
of $\mathsf{U}\mathfrak{g}$: 
$\ds \mathsf{U}\mathfrak{g}[[ t ]]
=\frac{\mathsf{T}\mathfrak{g}}{<xy-yx-t[x,y]>}$.
This is a Hopf algebra with $\Delta$, $\epsilon $ and $S$ as for 
$\mathsf{U}\mathfrak{g}$.
\item Let $\{X_i\ ;\ i=1,\ldots , n\}$ be a basis of $\mathfrak{g}$ 
and $\buildrel{X_i}_{}^{\rightarrow}$ (resp.
$\buildrel{X_i}_{}^{\leftarrow}$) be the left (resp. right) invariant 
vector fields on $G$ associated with $X_i$. For $\lambda\in [0,1]$
we consider the following actions of $B=\mathsf{U}\mathfrak{g}[[ t ]]$ 
on $C=\mathcal{C}^\infty (G)$:
\begin{enumerate}
\item $(X_i\rightharpoonup f) (x)=t
(\lambda-1)(\buildrel{X_i}_{}^\rightarrow\cdot f)(x)$
\item $(f\leftharpoonup X_i) (x)=t\lambda
(\buildrel{X_i}_{}^\leftarrow\cdot f) (x)$.
\end{enumerate}
\end{itemize}

\begin{lem}
These actions define on $\mathcal{C}^\infty (G)$ a $B$-bimodule 
algebra structure.
\end{lem}

\begin{dfn}
We denote by $\star_\lambda$ the L-R smash product on 
$\mathcal{C}^\infty (G)\otimes \Pol(\mathfrak{g}^\ast)$ given by this
$B$-bimodule algebra structure on $\mathcal{C}^\infty (G)$.
\end{dfn}

\begin{prop}
For $G=\mathbb{ R}^n$, $\star_{1/2}$ is the Moyal (Weyl ordered) 
star product, $\star_0$ is the standard ordered star product and in general 
$\star_\lambda$ is called $\lambda$-ordered star product on
$\mathbb{R}^{2n}$ \cite{Pflm99a}.
\end{prop}

\begin{rmk}
{\rm For a general Lie group $G$, $\star_\lambda$ gives in the generic case
new deformation quantization formulas on $T^*G$.
It would be interesting to study the properties of these $\star_\lambda$ 
for a noncommutative $G$ and their relations with the star products that 
are known. In particular $\star_{1/2}$ is formally different from the 
star product on $\mathcal{C}^\infty(T^\star G)$ given by S.~Gutt 
in \cite{Guts83a} but preliminary calculations seem to indicate that, in a
neighborhood of the unit of $G$, they are equivalent by a symplectomorphism.}
\end{rmk}

\subsection{Application to the quantization of symmetric spaces}

\begin{dfn}[\cite{Biep95a}] \label{SSSDEF} A {\em symplectic symmetric space} 
is a triple $(M, \omega,s)$, where $(M,\omega)$ is a smooth connected 
symplectic manifold and $s : M \times M \to M$ is a smooth map such that:

\begin{enumerate} 
\item[(i)] for all $x$ in $M$, the partial map $s_x : M \to M : y \mapsto 
s_x (y) := s(x,y)$ is an involutive symplectic diffeomorphism of $(M,\omega)$ 
called the {\em symmetry} at $x$. 
\item[(ii)] For all $x$ in $M$, $x$ is an isolated fixed point of $s_x$. 
\item[(iii)] For all $x$ and $y$ in $M$, one has $s_xs_ys_x=s_{s_x (y)}$. 
\end{enumerate}

Two symplectic symmetric spaces $(M,\omega,s)$ 
and $(M',\omega ',s')$ are {\em isomorphic} if there exists a symplectic 
diffeomorphism $\varphi: (M,\omega) \rightarrow (M',\omega')$ such that 
$\varphi s_x=s'_{\varphi (x)} \varphi$. 
\end{dfn}

\begin{dfn} Let $(\mathfrak{g},\sigma)$ be an {\em involutive algebra}, 
that is, $\mathfrak{g}$ is a finite dimensional real Lie algebra and 
$\sigma$ is an involutive automorphism of $\mathfrak{g}$. Let $\Omega$ be 
a skewsymmetric bilinear form on $\mathfrak{g}$. Then the triple 
$(\mathfrak{g},\sigma, \Omega)$ is called a {\em symplectic triple} 
if the following properties are satisfied: 

\begin{enumerate} 
\item Let $\mathfrak{g}=\mathfrak{k}\oplus \mathfrak{p}$ 
where $\mathfrak{k}$ (resp. 
$\mathfrak{p}$) is the $+1$ (resp. $-1$) eigenspace of $\sigma$. 
Then $[\mathfrak{p},\mathfrak{p}]=\mathfrak{k}$ 
and the representation of $\mathfrak{k}$ 
on $\mathfrak{p}$, given by the adjoint action, is faithful. 
\item $\Omega$ is a Chevalley 2-cocycle for the trivial representation 
of $\mathfrak{g}$ on $\mathbb{R}$ such that $\forall X \in \mathfrak{k}$, 
$i(X){\Omega}=0$. Moreover, the restriction of 
$\Omega$ to $\mathfrak{p} \times \mathfrak{p}$ is nondegenerate.
\end{enumerate} 

The dimension of $\mathfrak{p}$ defines the {\em dimension} of the triple. 
Two such triples $(\mathfrak{g}_i,\sigma_i,{\Omega}_i)$ $(i=1,2)$ are 
{\em isomorphic} if there exists a Lie algebra isomorphism 
$\psi :\mathfrak{g}_1\rightarrow\mathfrak{g}_2$ such that 
$\psi \circ \sigma_1 = \sigma_2 \circ \psi$ and $\psi^*{\Omega}_2={\Omega}_1$. 
\end{dfn}

\begin{prop}[\cite{Biep95a}] There is a bijective correspondence between 
the isomorphism classes of simply connected symplectic symmetric spaces 
$(M,\omega,s)$ and the isomorphism classes of symmetric triples 
$(\mathfrak{g},\sigma,{ \Omega})$. 
\end{prop}

\begin{dfn}
A symplectic symmetric space $(M,\omega,s)$ is called an {\em elementary 
solvable} symplectic symmetric space if its associated triple 
$(\mathfrak{g},\sigma,\Omega)$ is of the following type:
\begin{enumerate}
\item The Lie algebra $\mathfrak{g}$ is a split extension of Abelian Lie 
algebras $\mathfrak{a}$ and $\mathfrak{b}$~:
$$
\mathfrak{b}\longrightarrow\mathfrak{g}
\stackrel{\longleftarrow}{\longrightarrow}\mathfrak{a}.
$$
\item The automorphism $\sigma$ preserves the splitting 
$\mathfrak{g}=\mathfrak{b}\oplus\mathfrak{a}$.
\item There exists $\xi\in \mathfrak{k}^\ast$ such that $\Omega (X,Y)
=\delta \xi =<\xi,[X,Y]_\mathfrak{g}>$ (Chevalley 2-coboundary).
\end{enumerate}
\end{dfn}

For such an elementary solvable symplectic symmetric space there exists 
a global Darboux chart such that
$(M,\omega )\simeq (\mathfrak{p}=\mathfrak{l}\oplus \mathfrak{a},\Omega)$
\cite{Biep00a}.
So we have 
$$\mathcal{C}^\infty (M)\simeq\mathcal{C}^\infty (\mathfrak{p})
\simeq \mathcal{C}^\infty (\mathfrak{l})
\hat\otimes \mathcal{C}^\infty (\mathfrak{a})
\simeq \mathcal{C}^\infty (\mathfrak{l})
\hat\otimes \mathcal{C}^\infty (\mathfrak{l}^\ast)
\buildrel\supset_{a\simeq\mathfrak{l}^{\ast}}^{} 
\mathcal{C}^\infty (\mathfrak{l})\otimes \Pol(\mathfrak{l}^\ast)
\buildrel\simeq_{\mathfrak{l}\hbox{ abelian}}^{}
\mathcal{C}^\infty (\mathfrak{l})\otimes \mathsf{U}\mathfrak{l}$$

One can now define $\star_{1/2}$ (Moyal) on 
$\mathcal{C}^\infty (M)
\simeq \mathcal{C}^\infty (\mathfrak{l}\oplus \mathfrak{a})$ or,
using our preceding construction, on
$\mathcal{C}^\infty (\mathfrak{l})\otimes \mathsf{U}\mathfrak{l}.$

In order to have an {\em invariant} star product on $M$ under the action of $G$ 
(such that $\mathfrak{g}=\hbox{Lie}(G))$ P. Bieliavsky \cite{Biep00a} defines 
an integral transformation $S: \mathcal{C}^\infty (\mathfrak{l})
\rightarrow \mathcal{C}^\infty (\mathfrak{l})$ 
and then an invariant star product $\star_S$ by, for $T:=S\otimes \Id,$
$$(f\otimes a)\star_S (g\otimes b)
:=T^{-1}(T(f\otimes a)\star_{1/2} T(g\otimes b)).$$

Let us define $\quad f\bullet_S g:= S^{-1} (Sf . Sg),$
$\quad a\buildrel\rightharpoonup_{}^S f:=S^{-1}(a\rightharpoonup Sf)\quad$
and $\quad f\buildrel\leftharpoonup_{}^S a:= S^{-1} (Sf\leftharpoonup a).$

\begin{prop}[\cite{BiBo02}] 
$\star_S$ is the L-R smash product 
of $(\mathcal{C}^\infty (\mathfrak{l}),\bullet_S)$ by 
$\mathsf{U}\mathfrak{l}$ with the 
$\mathsf{U}\mathfrak{l}$-bimodule structure given by 
$\buildrel\rightharpoonup_{}^S$ 
and $\buildrel\leftharpoonup_{}^S$.
\end{prop}

\begin{rmk}
{\rm Since we were dealing with quantum groups in the first sections, 
we want to stress that the homogeneous (symmetric) spaces
involved here are strictly different from those appearing
in the quantum group approach of quantized homogeneous spaces
\cite{Driv93a}. Indeed, in the latter, the spaces come from
Poisson-Lie groups, so that the Poisson bracket has to be singular;
therefore this bracket (and a fortiori a star product deforming this 
bracket) cannot be invariant (otherwise it would be zero everywhere).
Here the Poisson brackets are invariant and regular.}
\end{rmk}

\noindent{\bf Acknowledgments.} 
This survey owes a lot to the insight shown by Mosh\'e
Flato in pushing forward the deformation quantization program, including
in its aspects related to quantum groups where the inputs of
Georges Pinczon and Murray Gerstenhaber were, as can be seen here,
very important. Thanks are also due to the referee for a number of
valuable comments.

\end{document}